\documentclass{amsart}
\usepackage{graphicx,subfigure}
\usepackage[english]{babel}
\usepackage[cp1251]{inputenc}
\usepackage{amssymb}
\usepackage{lscape}
\usepackage{amsthm}
\usepackage{inputenc}

\newtheorem{defn}{Definition}[section]
\newtheorem{thm}{Theorem}[section]

\newenvironment{dem}{\rm \trivlist \item[\hskip \labelsep{\it
      Proof}:]}{\par\nopagebreak \hfill $\Box$ \endtrivlist}

\def\ll{\mathcal{L}}

\begin{document}
\date{}
\author{L.M. Camacho, E.M. Ca\~{n}ete, J.R. G\'{o}mez, B.A. Omirov}

\title[p-filiform Leibniz algebras of maximum length]{\bf p-filiform Leibniz algebras of maximum length.}
\address{[L.M. Camacho --- E.M. Ca\~{n}ete --- J.R. G\'{o}mez] Dpto. Matem\'{a}tica Aplicada I.
Universidad de Sevilla. Avda. Reina Mercedes, s/n. 41012 Sevilla.
(Spain)} \email{lcamacho@us.es --- elisacamol@us.es --- jrgomez}
\address{[B.A. Omirov] Institute of Mathematics and Information Technologies
of Academy of Uzbekistan, 29, Do'rmon yo'li street., 100125, Tashkent (Uzbekistan)}
\email{omirovb@mail.ru}

\thanks{This work has been funded by Mathematics Institute and V Research Plan of Sevilla University, by the Grants of Junta de Andaluc\'{i}a, by
the Grants (RGA) No:11-018 RG/Math/$AS_{-}I$--UNESCO FR: 3240262715 and IMU/CDC-program.}%

\maketitle
\begin{abstract}
The descriptions (up to isomorphism) of naturally graded $p$-filiform Leibniz algebras and $p$-filiform ($p\leq 3$) Leibniz algebras of maximum length are known. In this paper we study the gradation of maximum length for $p$-filiform Leibniz algebras. The present work aims at the classification of complex $p$-filiform ($p \geq 4$) Leibniz algebras of maximum length.
\end{abstract}



\section{Introduction}
Leibniz algebras were introduced at the beginning of the 90s by J.-L. Loday in \cite{loday}, which are a ``non-commutative" generalization of Lie algebras.
The right multiplication operator on an element of a Leibniz
algebra is a derivation, which is a property inherited from Lie algebras.

Active investigation on Leibniz algebra theory shows that many results of the
theory of Lie algebras can be extended to Leibniz algebras. Distinctive properties of non-Lie
Leibniz algebras have also been studied \cite{Omirov1, Bar}.
For a Leibniz algebra there is a corresponding associated Lie algebra, which is the quotient algebra by the two-sided ideal generated by squares of elements of a Leibniz algebra (denoted by $I$).

From the theory of Lie algebras it is well known that the study of finite dimensional Lie algebras was
reduced to the nilpotent ones, due to Levi's theorem and Mal'cev's decomposition (see \cite{Jacobson, Mal'cev}). The case of Leibniz algebras is analogous to Levi's theorem
\cite{Bar}. Namely, a Leibniz algebra is decomposed into a semidirect sum of its solvable radical and a
semisimple Lie algebra.

The structure of solvable Lie algebra can be obtained from the structure of its nilradical \cite{Mub}. This approach has recently been extended to the case of Leibniz algebras \cite{Casas}. Therefore, the main problem when describing of finite-dimensional Leibniz algebras is the nilpotent radical. Thus, the study of nilpotent Leibniz algebras is a crucial problem.

Since the description of the set of $n$-dimensional nilpotent Leibniz algebras is an unsolvable task (even in the case of Lie algebras), we have to study nilpotent Leibniz algebras under certain conditions (conditions on index of nilpotency, various types of gradation, characteristic sequence etc.).

The well-known gradations of nilpotent Lie and Leibniz algebras are very
helpful when investigating of the properties of those algebras without
restrictions on the gradation. Indeed, we can always choose an homogeneous basis and
thus the gradation allows to obtain more explicit conditions for the structural constants.
Moreover, such gradation is useful for the investigation of cohomologies for the considered algebras, because it induces the corresponding gradation of the group of cohomologies.

The concept of length of a Lie algebra was introduced by G\'{o}mez, Jim\'{e}nez-Merch\'{a}n and Reyes in \cite{Jimenez1}, \cite{Gomez}. Where, they distinguished an interesting family: algebras admitting a gradation with the greatest possible number of non-zero subspaces. Actually, the gradations with a large number of non-zero subspaces enable us to describe the multiplication on the algebra more exactly. They called such algebras \emph{algebras of maximum length}.
In fact, they only consider the connected gradation. There exist non connected algebras with the
greatest possible number of non-zero subspaces. Nevertheless, according to G\'{o}mez et al. the
notion of algebras of maximum length has already been used.

In \cite{Jobir}, \cite{NGMio} - \cite{3-filiformLeib} the classification of $p$-filiform Leibniz algebras of maximum length for $0 \leq p \leq 3$ is already closed.

The present paper aims at classification of $n$-dimensional $p$-filiform Leibniz algebras of maximum length with $n$ and $p$ generic.

Throughout the paper we consider finite-dimensional vector spaces and non-split algebras
over the field of the complex numbers. Moreover, we have omitted null products in the multiplication table of an algebra.

\

\section{Preliminares}
 Recall \cite{loday} that an algebra $\ll$ over a field $F$ is called a Leibniz algebra if it satisfies the following Leibniz identity:
$$[x,[y,z]]=[[x,y],z]-[[x,z],y], \forall x,y,z \in \ll $$
where $[ - , - ]$ denotes the multiplication of an algebra $\ll$.

Let $\ll$ be a Leibniz algebra, then $\ll$ is naturally filtered by the descending central sequence $$\ll^1=\ll, \ \ll^{k+1}=[\ll^k,\ll], \ k\geq 1.$$

An nilpotent Lebniz algebra $\ll$ has {\it nilindex} equal to $s$ if $s$ is the minimum integer such that $\ll^s \neq \{0\}$ and $\ll^{s+1}=\{0\}.$


We denote by $R_x$ the operator of right multiplication on element $x$, i.e., $R_x: \ll \rightarrow \ll$ such that $R_x(y)=[y,x]$ for any $y \in \ll$.


 Let $x$ be an element of the set $\ll \setminus \ll^2$. For the nilpotent operator $R_x$ we define a descending sequence $C(x)=(n_1,n_2, \dots, n_k)$, which consists of the dimensions of the Jordan blocks of the operator $R_x$. In the set of such sequences we consider the lexicographic order, that is,
  $C(x)=(n_1,n_2, \dots, n_k)< C(y)=(m_1, m_2, \dots, m_s)$ if and only if there exists $i \in \mathbb{N}$ such that $n_j=m_j$ for some $j<i$ and $n_i<m_i$.

  \begin{defn}\label{def:char.seq}
  The sequence $C(\ll)=\max C(x)_{x \in \ll \setminus \ll^2}$ is called the characteristic sequence of the algebra $\ll$.
  \end{defn}

  Let $\ll$ be an $n$-dimensional nilpotent Leibniz algebra and $p$ be an non negative integer ($p<n$).
  \begin{defn}\label{def:p-fili}
  The Leibniz algebra $\ll$ is called $p$-filiform if $C(\ll)=(n-p,\underbrace{1,\dots,1}_{p})$.
  \end{defn}


A Leibniz algebra $\ll$ is $\mathbb{Z}$-graded if $\ll=\oplus_{i \in \mathbb{Z}}V_i$,
   where $[V_i,V_j]\subseteq V_{i+j}$ for some $i,j \in \mathbb{Z}$ with a finite number of non null spaces $V_i$.

  We will say that a $\mathbb{Z}$-graded nilpotent Leibniz algebra $\ll$ admits a \emph{connected gradation} if $\ll=V_{k_1}\oplus V_{k_1+1} \oplus \dots \oplus V_{k_1+t}$ and $V_{k_1+i}\neq <0>$ for some $i$ $(0 \leq i \leq t)$.
  \begin{defn}\label{def:length}
  The number $l( \oplus \ll)=l(V_{k_1}\oplus V_{k_1+1} \oplus \dots \oplus V_{k_1+t})=t+1$ is called the
length of the gradation, where $ \oplus \ll$ is a connected gradation. The gradation $\oplus \ll$ has maximum length if $l(\oplus \ll)=dim( \ll)$.
  \end{defn}

  We define the length of an algebra $\ll$ as follows
  \begin{center}
 $l(\ll)=\max \{l(\oplus \ll)  \hbox{ such that } \oplus \ll = V_{k_1}\oplus \dots \oplus V_{k_t}\hbox{ is a connected gradation}\}.$
 \end{center}

An algebra $\ll$ is called {\it of maximum length} if $l(\ll)=dim (\ll)$.

Thus, we resume the properties of gradation of maximum length:
   \begin{equation}\label{grad}
   \begin{cases}
   [V_i,V_j] \subseteq V_{i+j},\\
   \hbox{the subspaces of the gradation are not empty,}\\
   \hbox{the dimension of each subspace of the gradation equals }1,\\
   \hbox{all subindices are different from each other.}
   \end{cases}
   \end{equation}

We define another type of gradation below.

Given an $n$-dimensional Leibniz algebra $\ll$ with nilindex $s$, put
$\ll_i=\ll^i/\ll^{i+1}$ with $ 1\leq i\leq s$ and $gr \ll =\ll_1 \oplus \ll_2
\oplus \cdots \oplus \ll_{s}.$ Then $[\ll_i,\ll_j]\subseteq \ll_{i+j}$
and we obtain the graded algebra $gr \ll$. If $gr \ll$ and $\ll$ are
isomorphic, $gr \ll\cong \ll,$ we say that $\ll$ is naturally graded.

The classification of naturally graded $p$-filiform Lie algebras was done by Cabezas y Pastor in \cite{Pastor}, which is given below.
\begin{thm}\label{p-filNG}
 Let $\ll$ be a $n$-dimensional naturally graded $p$-filiform Lie algebra, with $p>1,$ $n \geq \max\{3p-1,p+8\}$ and $3 \leq r_1 < r_2 < \dots <r_{p-1} \leq n-p$ odds. Therefore:
\small\begin{itemize}
\item{If $r_{p-1}=n-p,$ then $\ll$ is isomorphic to $L(n,r_1,r_2,\dots,r_{p-2},n-p)$,}
\item{If $r_{p-1}=n-p-1,$ then $\ll$ is isomorphic to $L(n,r_1,r_2,\dots,r_{p-2},n-p-1)$ or $\tau(n,r_1,r_2,\dots,r_{p-2},n-p-1)$,}
\item{If $r_{p-1}=n-p-2,$ then $\ll$ is isomorphic to $L(n,r_1,r_2,\dots,r_{p-2},n-p-2),$ $Q(n,r_1,r_2,\dots,r_{p-2},n-p-2)$ or  $\tau(n,r_1,r_2,\dots,r_{p-2},n-p-2)$,}
\item{If $2p-1 \leq r_{p-1} \leq n-p-3,$} then
\begin{itemize}
\item{if $n-p$ is odd,} then $\ll$ is isomorphic to $L(n,r_1,r_2,\dots,r_{p-1})$ or $Q(n,r_1,r_2,\dots,r_{p-1})$,
\item{if $n-p$ is even,} then $\ll$ is isomorphic to $L(n,r_1,r_2,\dots,r_{p-1}).$
\end{itemize}
\end{itemize}
 where
$$\small\begin{array}{l}
 L(n,r_1,r_2,\dots,r_{p-1})\ :\\
 \begin{cases}
                [x_{0},x_{i}]=x_{i+1}, & 1\leq i \leq n-p-1, \\
                [x_{i},x_{r_j-i}]=(-1)^{i-1}y_{j}, & 1\leq i \leq \frac{r_j-1}{2}, \ 1 \leq j \leq p-1.\\
          \end{cases}

 \\[10mm]
 Q(n,r_1,r_2, \dots, r_{p-1})\ :\\
 \begin{cases}
                [x_{0},x_{i}]=x_{i+1}, & 1\leq i \leq n-p-1, \\
                [x_{i},x_{r_j-i}]=(-1)^{i-1}y_{j}, & 1\leq i \leq \frac{r_j-1}{2}, \ 1 \leq j \leq p-1,\\
                [x_{i},x_{n-p-i}]=(-1)^{i-1}x_{n-p}, & 1\leq i \leq \frac{n-p-1}{2}.
          \end{cases}
 \\[10mm]
\tau(n,r_1,r_2, \dots,r_{p-2},n-p-1) \ : \\
\begin{cases}
                [x_{0},x_{i}]=x_{i+1}, & 1\leq i \leq n-p-1, \\
                [x_{i},x_{r_j-i}]=(-1)^{i-1}y_{j},  & 1\leq i \leq \frac{r_j-1}{2}, \ 1 \leq j \leq p-2,\\
                [x_{i},x_{n-p-1-i}]=(-1)^{i-1}(x_{n-p-1}+y_{p-1}), & 1\leq i \leq \frac{n-p-2}{2},\\
                [x_{i},x_{n-p-i}]=(-1)^{i-1}\frac{n-2i-p}{2}x_{n-p}, & 1\leq i \leq \frac{n-p-2}{2},\\
                [x_1,y_{p-1}]=\frac{p+2-n}{2}x_{n-p}.
          \end{cases}
\\[20mm]
\tau(n,r_1,r_2, \dots,r_{p-2},n-p-2) \ : \\
\begin{cases}
                [x_{0},x_{i}]=x_{i+1}, & 1\leq i \leq n-p-1, \\
                [x_{i},x_{r_j-i}]=(-1)^{i-1}y_{j},  & 1\leq i \leq \frac{r_j-1}{2}, \ 1 \leq j \leq p-2,\\
                [x_{i},x_{n-p-2-i}]=(-1)^{i-1}(x_{n-p-2}+y_{p-1}), & 1\leq i \leq \frac{n-p-3}{2},\\
                [x_{i},x_{n-p-1-i}]=(-1)^{i-1}\frac{n-p-1-2i}{2}x_{n-p-1}, & 1\leq i \leq \frac{n-p-3}{2},\\
                [x_{i},x_{n-p-i}]=(-1)^{i}(i-1)\frac{n-p-1-i}{2}x_{n-p}, & 2\leq i \leq \frac{n-p-1}{2},\\
                [x_i,y_{p-1}]=\frac{p+3-n}{2}x_{n-p-2+i}, & 1 \leq i \leq 2.
          \end{cases}
\end{array}$$
with $\{x_0,x_1,\dots,x_{n-p},y_1,\dots,y_{p-1}\}$ a basis.
\end{thm}

For naturally graded $p$-filiform non-Lie Leibniz algebras the result obtained is the following.

\begin{thm} \cite{NGp-F} \label{thm22}
Let $\ll$ be a $n$-dimensional $p$-filiform non-Lie Leibniz algebra, with $n-p \geq 4.$ Then $\ll$ is isomorphic to one of the following algebras
\begin{itemize}
\item[]{If $p$ is even:}
$$M^1:= \begin{cases}
[e_i,e_1]=e_{i+1}, \quad 1 \leq i \leq n-p-1,\\
[e_1,f_j]=f_{\frac{p}{2}+j}, \quad 1 \leq j \leq \frac{p}{2}.
\end{cases}$$

$$M^2:= \begin{cases}
[e_i,e_1]=e_{i+1}, \quad 1 \leq i \leq n-p-1,\\
[e_1,f_1]=e_2+f_{\frac{p}{2}+1},\\
[e_i,f_1]=e_{i+1}, \quad 1 \leq i \leq n-p-1,\\
[e_1,f_j]=f_{\frac{p}{2}+j}, \quad 2 \leq j \leq \frac{p}{2}.\\

\end{cases}$$

\item[]{If $p$ is odd:}
$$M^3:= \begin{cases}
[e_i,e_1]=e_{i+1}, \quad 1 \leq i \leq n-p-1,\\
[e_1,f_j]=f_{\lfloor\frac{p}{2}\rfloor+j}, \quad 1 \leq j \leq \lfloor \frac{p}{2}\rfloor,\\
[e_i,f_{\lfloor \frac{p}{2} \rfloor}+1]=e_{i+1}, \quad 1 \leq i \leq n-p-1.
\end{cases}$$
\end{itemize}
with $\{e_1,\dots,e_{n-p},f_1,\dots,f_{p}\}$ a basis.
\end{thm}

\section{$p$-filiform Leibniz algebras of maximum length}

\

In order to achieve our goal we use the following algorithm:

1. Firstly, we extend the naturally graded $p$-filiform Leibniz algebras by using the natural gradations. In this way, we can distinguish two cases: the natural graded $p$-filiform Lie algebras and the natural graded $p$-filiform non-Lie Leibniz algebras.

2. After that, we construct an homogeneous basis of a graded algebra of maximum length with respect to the basis of naturally gradation.

3. Finally, we classify the $p$-filiform Leibniz algebra in an homogeneous basis of maximum length.

\subsection{Extension of Lie algebras}

\

\

In this subsection we prove that there is no $n$-dimensional $p$-filiform Leibniz algebra of maximum length in the non split case, for $p \geq 4$ and $n \geq \max\{3p-1,p+8\}.$ Note that the study of the particular case for $n < \max\{3p-1,p+8\}$ can be found in \cite{mitesis}.

Further we will denote by $\widetilde{\ll}$ the extension of $\ll,$ that is, the family of algebras with the table of multiplication whose naturally graded algebra is $\ll.$

\begin{thm} \label{thm31}
 Let $\ll$ be a $n$-dimensional $p$-filiform Leibniz algebra, whose naturally graded associated algebra is isomorphic to $L(n,r_1,r_2,\dots,r_{p-1})$ or  $Q(n,r_1,r_2, \dots, r_{p-1}),$ with $p \geq 4,$ $n \geq \max\{3p-1,p+8\}$ and $3 \leq r_1 < r_2 < \dots <r_{p-1} \leq n-p.$ Then the algebra $\ll$ does not admit a gradation of maximum length.
\end{thm}
\begin{dem}

Note that, we have $p \geq 4$, $n \geq 11$ and $3 \leq r_1 <r_2 \dots < r_{p-1}\leq n-p,$ where all $r_i$ are odds.

Let us suppose that $\ll$ admits a gradation of maximum length.

\

 $\fbox{\emph{Study of the extension $\widetilde{L}(n,r_1,r_2,\dots,r_{p-1}).$}}$

\

It is easy to see that natural gradation of the algebra $L(n,r_1,r_2,\dots,r_{p-1}$ is $L_1=\langle x_0,x_1\rangle,$ $L_i=\langle x_i\rangle$ with $2 \leq i \leq n-p$ and $i \neq r_j$ and $L_{r_j}=\langle x_{r_j},y_j\rangle$ with $1 \leq j \leq p-1.$ Therefore the law of $\widetilde{L}(n,r_1,r_2,\dots,r_{p-1})$ is defined by the following products, where the asterisks (*) denote the corresponding structural constants:
    $$\small\begin{cases}
  [x_0,x_0]=(*)x_3+\dots+(*)x_{n-p}+(*)y_1+\dots+(*)y_{p-1},\\
  [x_0,x_i]=x_{i+1}+(*)x_{i+2}+ \dots +(*)x_{n-p}+(*)y_1+\dots+(*)y_{p-1}, & 1 \leq i \leq r_1-2,\\
  [x_0,x_i]=x_{i+1}+(*)x_{i+2}+ \dots +(*)x_{n-p}+(*)y_2+\dots+(*)y_{p-1}, & r_1-1 \leq i \leq r_2-2,\\
   \hspace{5 cm} \vdots  & \hspace{1 cm} \vdots\\
  [x_0,x_i]=x_{i+1}+(*)x_{i+2}+ \dots +(*)x_{n-p}+(*)y_{p-1}, & r_{p-2}-1 \leq i \leq r_{p-1}-2,\\
  [x_0,x_i]=x_{i+1}+(*)x_{i+2}+ \dots +(*)x_{n-p}, & r_{p-1}-1 \leq i \leq n-p-1,\\
  [x_i,x_0]=-x_{i+1}+(*)x_{i+2}+ \dots +(*)x_{n-p}+(*)y_1+\dots+(*)y_{p-1}, & 1 \leq i \leq r_1-2 ,\\
  [x_i,x_0]=-x_{i+1}+(*)x_{i+2}+ \dots +(*)x_{n-p}+(*)y_2+\dots+(*)y_{p-1}, & r_1-1 \leq i \leq r_2-2,\\
  \hspace{5 cm} \vdots  & \hspace{1 cm} \vdots\\
  [x_i,x_0]=-x_{i+1}+(*)x_{i+2}+ \dots +(*)x_{n-p}+(*)y_{p-1}, & r_{p-2}-1 \leq i \leq r_{p-1}-2,\\
  [x_i,x_0]=-x_{i+1}+(*)x_{i+2}+ \dots +(*)x_{n-p}, & r_{p-1}-1 \leq i \leq n-p-1,\\
  [x_i,x_{r_j-i}]=(*)x_{r_j+1}+\dots+(*)x_{n-p}+(-1)^{i-1}y_j+(*)y_{j+1}+\dots+(*)y_{p-1}, & 1 \leq i \leq \frac{r_j-1}{2}, \\
  & 1 \leq j \leq p-1, \ r_j+1 \geq 4,\\
  [x_{r_j-i},x_i]=(*)x_{r_j+1}+\dots+(*)x_{n-p}+(-1)^{i}y_j+(*)y_{j+1}+\dots+(*)y_{p-1}, & 1 \leq i \leq \frac{r_j-1}{2},\\
  & 1 \leq j \leq p-1, \ r_j+1 \geq 4,\\
  [x_{i},x_{j}]=(*)x_{i+j+1}+\dots(*)x_{n-p}+(*)y_{r_k}+\dots+(*)y_{p-1}, & 1\leq i+j \leq r_k, \\
   & 1 \leq k \leq p-1,\\
    [x_{i},y_{j}]=(*)x_{i+r_j+1}+\dots+(*)x_{n-p}+(*)y_{r_k}+\dots+(*)y_{p-1}, & 0 \leq i \leq n-p-1, \\
   & 1 \leq j \leq p-1\\
  [y_i,y_j]=(*)x_{r_i+r_j+1}+\dots+(*)x_{n-p}+(*)y_{r_k}+\dots+(*)y_{p-1}, & 1 \leq i,j \leq p-4,\\
  & 7 \leq r_k \leq p-1, \\
  & r_i+r_j+1 \leq r_k \leq p-1.
  \end{cases}$$
Without loss of generality, one can assume that the general form of homogeneous generators of the gradation of maximum length are
\begin{equation}\label{generadores}
\widetilde{x_s}=x_0+\sum_{i=1}^{n-p}a_ix_i+\sum_{j=1}^{p-1}b_jy_j \quad \hbox{ and } \quad \widetilde{x_t}=A_0x_0+x_1+\sum_{i=2}^{n-p}A_ix_i+\sum_{j=1}^{p-1}B_jy_j,
\end{equation}
where $det\left( \begin{array}{cc}
1 & a_1\\
A_0 & 1
\end{array} \right) \neq 0.$

Below we present the straightforward consequences of the generated elements of $\widetilde{L}(n,r_1,r_2,\dots,r_{p-1})$ via the above two vectors:
{\small\begin{align*}
[\widetilde{x_s},\widetilde{x_s}]&=(*)x_3+\dots+(*)x_{n-p}+(*)y_1+\dots+(*)y_{p-1},\\
[\widetilde{x_t},\widetilde{x_t}]&=(*)x_3+\dots+(*)x_{n-p}+(*)y_1+\dots+(*)y_{p-1},\\
[\widetilde{x_s},\widetilde{x_t}]&=(1-a_1A_0)x_2+(*)x_3+\dots+(*)x_{n-p}+(*)y_1+\dots+(*)y_{p-1},\\
[[[\widetilde{x_s},\widetilde{x_t}],\underbrace{\widetilde{x_s}], \dots \widetilde{x_s}}_{i-times}]&=(-1)^i(1-a_1A_0)x_{i+2}+(*)x_{i+3}+\dots+(*)x_{n-p}+(*)y_1+\dots+(*)y_{p-1}\\
& \hbox{ for } 1 \leq i \leq r_1-3,\\
[[[\widetilde{x_s},\widetilde{x_t}],\underbrace{\widetilde{x_s}], \dots \widetilde{x_s}}_{i-times}]&=(-1)^i(1-a_1A_0)x_{i+2}+(*)x_{i+3}+\dots+(*)x_{n-p}+(*)y_2+\dots+(*)y_{p-1}\\
& \hbox{ for } r_1-1 \leq i \leq r_2-3,\\
[[[\widetilde{x_s},\widetilde{x_t}],\underbrace{\widetilde{x_s}], \dots \widetilde{x_s}}_{i-times}]&=(-1)^i(1-a_1A_0)x_{i+2}+(*)x_{i+3}+\dots+(*)x_{n-p}+(*)y_j+\dots+(*)y_{p-1} \\
& \hbox{ for } r_{j-1}-1 \leq i \leq r_j-3, \ 2 \leq j \leq p-1,\\
[[[\widetilde{x_s},\widetilde{x_t}],\underbrace{\widetilde{x_s}], \dots \widetilde{x_s}}_{i-times}]&=(-1)^i(1-a_1A_0)x_{i+2}+(*)x_{i+3}+\dots+(*)x_{n-p}\\
 &\hbox{ for } r_{p-1}-1 \leq i \leq n-p-2.\\
[[[\widetilde{x_s},\widetilde{x_t}],\underbrace{\widetilde{x_s}], \dots \widetilde{x_s}}_{(r_k-2)-times}]&=(-1)^{r_k-2}(1-a_1A_0)x_{r_k}+(*)x_{i+3}+\dots+(*)x_{n-p}+(-1)^{r_k-2}(1-a_1A_0)a_1y_k+\\
& +(*)y_{k+1}+\dots+(*)y_{p-1} \hbox{ for } 1 \leq k \leq p-1.
\end{align*}}

Let us take the new homogeneous basis constructed by the following vectors:
$$\begin{array}{ll}
z_0=\widetilde{x_s}, & z_1=\widetilde{x_t},\\
z_2=[z_0,z_1], & \\
z_i=[z_{i-1},z_0], & 3 \leq i \leq n-p, \\
p_j=[z_1,z_{r_j-1}], & 1 \leq j \leq p-1,
\end{array}$$
where
$$\begin{array}{lll}
p_j=[z_1,z_{r_j-1}]&=&(-1)^{r_j-3}(1-a_1A_0)A_0x_{r_j}+(*)x_{r_j+1}+\dots+(*)x_{n-p}+\\
&+&(-1)^{r_1-3}(1-a_1A_0)y_j+(*)y_{j+1}+(*)y_{p-1}.
\end{array}$$

The matrix of the change of basis is
\begin{center}
\begin{landscape}
 $$\Large\left(\begin{array}{ccccccccccccccccc}
 1 & a_1 & a_2 & a_3 & \dots & a_{r_1} & \dots & a_{r_2} & \dots & a_{r_{p-1}} & \dots & a_{n-p} & b_1 & b_2 & \dots & b_{p-1}\\[2mm]
 A_0 & 1 & A_2 & A_3 & \dots & A_{r_1} & \dots & A_{r_2} & \dots & A_{r_{p-1}} & \dots & A_{n-p} & B_1 & B_2 & \dots & B_{p-1}\\[2mm]
 0 & 0& C_2 & (*) & \dots & (*) & \dots & (*) & \dots & (*) & \dots & (*) & (*) & (*) & \dots & (*)\\[2mm]
  0 & 0& 0 & C_3 & \dots & (*) & \dots & (*) & \dots & (*) & \dots & (*) & (*) & (*) & \dots & (*)\\[2mm]
 \vdots &  \vdots &  \vdots &  \vdots &  \vdots &  \vdots &  \vdots &  \vdots &  \vdots &  \vdots &  \vdots &  \vdots &  \vdots &  \vdots &  \vdots &  \vdots \\[2mm]
  0 & 0& 0 & 0 & \dots & C_{r_1} & \dots & (*) & \dots & (*) & \dots & (*) & C_{r_1}a_1 & (*) & \dots & (*)\\[2mm]
   \vdots &  \vdots &  \vdots &  \vdots &  \vdots &  \vdots &  \vdots &  \vdots &  \vdots &  \vdots &  \vdots &  \vdots &  \vdots &  \vdots &  \vdots &  \vdots\\[2mm]
  0 & 0& 0 & 0 & \dots & 0 & \dots & C_{r_2} & \dots & (*) & \dots & (*) & 0 & C_{r_2}a_1 & \dots & (*)\\[2mm]
   \vdots &  \vdots &  \vdots &  \vdots &  \vdots &  \vdots &  \vdots &  \vdots &  \vdots &  \vdots &  \vdots &  \vdots &  \vdots &  \vdots &  \vdots &  \vdots\\[2mm]
  0 & 0& 0 & 0 & \dots & 0 & \dots & 0 & \dots & C_{r_{p-1}} & \dots & (*) & 0 & 0 & \dots & C_{r_{p-1}}\\[2mm]
  \vdots &  \vdots &  \vdots &  \vdots &  \vdots &  \vdots &  \vdots &  \vdots &  \vdots &  \vdots &  \vdots &  \vdots &  \vdots &  \vdots &  \vdots &  \vdots \\[2mm]
   0 & 0& 0 & 0 & \dots & 0 & \dots & 0 & \dots & 0 & \dots & C_{n-p} & 0 & 0 & \dots & 0\\[2mm]
   0 & 0& 0 & 0 & \dots & -C_{r_1}A_0 & \dots & (*) & \dots & (*) & \dots & (*) & -C_{r_1} & (*) & \dots & (*)\\[2mm]
   0 & 0& 0 & 0 & \dots & 0 & \dots & -C_{r_2}A_0 & \dots & (*) & \dots & (*) & 0 & -C_{r_2} & \dots & (*)\\[2mm]
  \vdots &  \vdots &  \vdots &  \vdots &  \vdots &  \vdots &  \vdots &  \vdots &  \vdots &  \vdots &  \vdots &  \vdots &  \vdots &  \vdots &  \vdots &  \vdots \\[2mm]
   0 & 0& 0 & 0 & \dots & 0 & \dots & 0 & \dots & -C_{r_{p-1}}A_0 & \dots & (*) & 0 & 0 & \dots & -C_{r_{p-1}}
 \end{array}\right)$$
 \end{landscape}
 \end{center}
  where $C_i=(-1)^i(1-a_1A_0),$ for $2 \leq i \leq r_p-1.$

 Note that this matrix has rank equal to $n$ (because of $1-a_1A_0\neq 0$).

We put $z_0 \in V_{k_s}$ and $z_1 \in V_{k_t}$. Then, according to the definition of gradation of maximum length, we derive:
$$\ll=V_{k_s} \oplus V_{k_t} \oplus V_{k_t+k_s} \oplus \dots \oplus V_{k_t+(n-p-1)k_s} \oplus V_{2k_t+(r_1-2)k_s} \oplus \dots \oplus V_{2k_t+(r_{p-1}-2)k_s}.$$ This gradation is connected if and only if $k_s=\pm 1.$ Since the cases $k_s=1$ and $k_s=-1$ are equivalent, one can assume $k_s=1.$

To make the reasoning simple, we sometimes shall continue the notation $k_s,$ even though that $k_s=1.$

    Our next objective is to analyze the value of $k_t.$ We distinguish the following cases:
\begin{itemize}
     \item{Let $k_t>0.$} Then the properties (\ref{grad}) are satisfied if and only if $k_t=2.$ In this case we conclude $k_t+k_s< 2k_t+(r_1-2)k_s <k_t+(n-p-1)k_s,$ that is, $V_{2k_t+(r_1-2)k_s}=\langle p_1,z_m \rangle \hbox{ with } 2\leq m \leq n-p,$ giving rise to a contradiction with the assumption of maximum length.

\item{Let $k_t<0.$} Then we consider the following subcases:
\begin{itemize}
\item{If $k_t=-n+p+1,$} then $2k_t+(r_i-2)k_s \leq 2(-n+p+1)+n-p-2=-n+p \leq 0$ for $1 \leq i \leq p-1.$ Hence we can affirm that all subspaces $V_{2k_t+(r_i-2)k_s}$ have negative subindices. Therefore, due to the connectedness of the gradation, we obtain:
$$\begin{array}{lll}
distance(p_i,p_{i+1})&=&distance(2k_t+r_i-2,2k_t+r_{i+1}-2)=1,\\
&\Rightarrow &  distance(r_i,r_{i+1})=1,
\end{array}$$
which is impossible since the parameters are odd $r_i$ for $1 \leq i \leq p-1.$
\item{If $k_t>-n+p+1,$} then there exists an $z_t,$ with $1 \leq t \leq n-p,$ such that $z_t \in V_1=\langle z_0 \rangle.$ However, it is impossible because of $z_t$ is generator. Hence we get a contradiction.
\item{If $k_t<-n+p+1,$} then it is easy to see that $V_0=\langle 0 \rangle,$ which contradicts the properties (\ref{grad}).
\end{itemize}
\end{itemize}
Thus, it has been proved that there is no algebra of  maximum length among the extension of the family $L(n,r_1,r_2,\dots,r_{p-1}).$

Furthermore, the above arguments can be used for some admissible value of $r_{p-1},$ that is, this proof includes the particular cases $\widetilde{L}(n,r_1,r_2,\dots,r_{p-2},n-p),$ $\widetilde{L}(n,r_1,r_2,\dots,r_{p-2},n-p-1)$ and $\widetilde{L}(n,r_1,r_2,\dots,r_{p-2},n-p-2).$

\

$\fbox{\emph{Study of the extension $\widetilde{Q}(n,r_1,r_2,\dots,r_{p-1}).$}}$

 \

This case is analogous to case $\widetilde{L}(n,r_1,r_2,\dots,r_{p-1}).$ The difference is only in the construction of the homogeneous basis of gradation of maximum length.

Let us consider the following cases:

\

\fbox{\textbf{Case 1: $1+a_1 \neq 0.$}} Then we take as a basis the following vectors:
\begin{align*}
z_0&=\widetilde{x_s},\\
z_1&=\widetilde{x_t},\\
z_2&=[z_0,z_1],\\
z_i&=[z_{i-1},z_0], \hbox{ for } 3 \leq i \leq n-p,\\
p_j&=[z_1,z_{r_j-1}], \hbox{ for } 1 \leq j \leq p-1,
\end{align*}
where
$$
p_j=A_0(-1)^{r_j-3}(1-a_1A_0)x_{r_j}+(*)x_{r_j+1}+\dots+(*)x_{n-p}+$$
$$(-1)^{r_1-3}(1-a_1A_0)y_j +(*)y_{j+1}+\\
 +\dots+(*)y_{p-1}.
$$

Since the chosen basis is the same as in the study of the extension $\widetilde{L}(n,r_1,r_2,\dots,r_{p-1}),$
for the same reasons we conclude that there is not any algebra of maximum length in this case.

\

\fbox{\textbf{Case 2: $1+a_1=0.$}} It is clear that $A_0 \neq -1$ (because of $1-a_1A_0 \neq 0$). Therefore, we can take the new basis defined by the following vectors:
 \begin{align*}
z_0&=\widetilde{x_s},\\
z_1&=\widetilde{x_t},\\
z_2&=[z_0,z_1],\\
z_i&=[z_{i-1},z_0], \hbox{ for } 3 \leq i \leq n-p-1,\\
z_{n-p}&=[z_{n-p-1},z_1]=(-1)^{n-p-2}(1-a_1A_0)(1+A_0)x_{n-p},\\
p_j&=[z_1,z_{r_j-1}], \hbox{ for } 1 \leq j \leq p-1.
\end{align*}
It is not difficult to check that the matrix of basis transformation is non-singular.

The associated gradation of maximum length is
$$V_{k_s} \oplus V_{k_t} \oplus V_{k_t+k_s} \oplus \dots \oplus V_{k_t+(n-p-2)k_s} \oplus V_{2k_t+(n-p-2)k_s}\oplus V_{2k_t+(r_1-2)k_s} \oplus \dots \oplus V_{2k_t+(r_{p-1}-2)k_s}.$$

Note that this gradation also satisfies the pro\-per\-ties (\ref{grad}). As stated above, without loss of generality, we can assume $k_s=1.$
We consider the following subcases:

    \begin{itemize}
    \item{Let $k_t>0.$} Then, by considering the properties (\ref{grad}), we conclude $k_t=2$. So, we have $k_s=1,$ $k_t=2$ and $3 \leq r_1 < n-p-2,$ that is, $5 \leq 2k_t+(r_1-2)k_s < n-p.$ Thus, one can assert the existence of $z_t$ with $4 \leq t \leq n-p-1,$ such that  $V_{2k_t+(r_1-2)k_s}=\langle p_1, z_t \rangle,$ which contradicts the assumption of maximum length (see the properties (\ref{grad})).

    \item{Let $k_t < 0$.} Then from the properties of gradation of maximum length we get $k_t=-n+p+2.$ Moreover, since $r_1$ and $r_2$ are odd we conclude  $$distance(p_1,p_2)=distance(2k_t+(r_1-2)k_s,2k_t+(r_2-2)k_s)=r_2-r_1 > 1.$$ Therefore, we obtain a contradiction with the assumption of maximum length again.
    \end{itemize}
\end{dem}
%

The next theorem is proved by applying the same methods and arguments as in the proof of Theorem \ref{thm31}.

\begin{thm}
Let $\ll$ be a $n$-dimensional $p$-filiform Leibniz algebra, whose naturally gradation leads to an algebra isomorphic to $\tau(n,r_1,r_2,\dots,n-p-1)$ or  $\tau(n,r_1,r_2, \dots, n-p-2),$ with $p \geq 4$ and $n \geq \max\{3p-1,p+8\}$. Then $\ll$ does not admit a gradation of maximum length.
\end{thm}

\

\subsection{Extension of Leibniz algebras}

\

\



In this subsection we study the description of $p$-filiform Leibniz algebras of maximum length from the extensions of naturally graded $p$-filiform non-Lie Leibniz algebras. The classification of Theorem \ref{thm22} leads to considerate the extensions of the algebras $M_1 - M_3$.

Firstly, we analyze the extension of the algebras $M^1$ and $M^2.$

\begin{thm}
Let $\ll$ be a $n$-dimensional $p$-filiform Leibniz algebra, with $n-p \geq 4,$ $p \geq 4$ and $p$ even. Then, the algebra $\ll$ is isomorphic to the one of the following pairwise non-isomorphic algebras:
$$\begin{array}{cc}
M^4(\alpha):\begin{cases}
[x_i,x_1]=x_{i+1},&  1 \leq i \leq n-p-1,\\
[x_1,y_i]=z_i, & 1 \leq i \leq \frac{p}{2},\\
[z_1,y_2]=[z_2,y_1]=\alpha x_{n-p},  & \alpha \in \{0,1\}.
\end{cases} &
M^5:\begin{cases}
[x_i,x_1]=x_{i+1}, & 1 \leq i \leq n-p-1,\\
[x_1,y_i]=z_i, & 1 \leq i \leq \frac{p}{2},\\
[y_1,y_2]=x_{n-p}.
\end{cases}
\end{array}$$
where $n$ is even and $\frac {n} {n-p} \in \mathbb{N}$ in the algebra $M^4(1).$
\end{thm}
\begin{dem} Similarly as above the generators of maximum length gradation of the algebra $\ll$ have the form:

\begin{align*}
\widetilde{x_s}&=e_1+ \sum_{i=2}^{n-p} a_ie_i+\sum_{i=1}^p b_jf_j,\\
y_j&=f_j+\sum_{k=1}^{n-p}c_{kj}e_k+\sum_{k=1, k\neq j}^{\frac{p}{2}}d_{kj}f_k \quad \hbox{ for } 1 \leq j \leq \frac{p}{2}.
\end{align*}

 $\fbox{\emph{Study of the extension $\widetilde{M^1}.$}}$

 \

 \

By considering the law of the algebra $M^1$ and the natural gradation we have: $$L_1=\langle e_1,f_1, \dots, f_{\frac{p}{2}} \rangle \oplus L_2=\langle e_2, f_{\frac{p}{2}+1}, \dots, f_p \rangle \oplus L_3=\langle e_3 \rangle \oplus \dots \oplus L_{n-p}=\langle e_{n-p} \rangle.$$

We construct the following new adapted basis:
$$x_1=\widetilde{x_s},\quad x_i=[x_{i-1},x_1], \ 2 \leq i \leq n-p, \quad
y_i, \ 1 \leq i \leq \frac{p}{2}, \quad z_i=[x_1,y_i], \ 1 \leq i \leq \frac{p}{2},
$$
 whose associated gradation is $$V_{k_s} \oplus V_{2k_s} \oplus \dots \oplus V_{(n-p)k_s} \oplus V_{k_1} \oplus V_{k_2} \oplus \dots \oplus V_{k_{\frac{p}{2}}} \oplus V_{k_1+k_s} \oplus V_{k_2+k_s} \oplus \dots \oplus V_{k_{\frac{p}{2}}+k_s}.$$

 Let us assume that this gradation has maximum length.

 We consider the products $[z_i,x_1]$ with $1 \leq i \leq n-p-1.$ Due to the law of the algebra $M_1$ we obtain  $[z_i,x_1]=c_{1i}e_3+(*)e_4+\dots+(*)e_{n-p},$ i.e., we can assume that $[z_i,x_1]=c_{1i}x_3, \ 1 \leq i \leq n-p-1.$

 On the other hand, by using the properties of the gradation we derive
 $$\begin{cases}
 [z_i,x_1] \in V_{2k_s+k_i}, \ \hbox{ for } 1 \leq i \leq n-p-1,\\
 x_3 \in V_{3k_s},
 k_s \neq k_i, \hbox{ for } 1 \leq i \leq n-p-1.
 \end{cases}$$

Therefor, $c_{1i}=0, \ 1 \leq i \leq n-p-1.$

By induction on a fixed $i$ and any $j$ and using the Leibniz identity, one can prove that $$[x_i,x_j]=0, \ 3 \leq i,j \leq n-p.$$

Let us analyze the products $[x_i,y_j], \ [y_j,x_i]$ for $1 \leq i \leq n-p,$ $1 \leq j \leq \frac{p}{2}$ and $(i,j)\neq (1,j).$

According to the law of the algebra $M_1$ we get
$[x_2,y_j]=(*)e_4+ \dots + (*)e_{n-p},$ that is, $[x_2,y_j]=Ax_m$ with $1 \leq j \leq \frac{p}{2}, \  4 \leq m \leq n-p$ and some coefficient $A$.

On the other hand, the properties of the gradation of maximum length deduce
$$\begin{cases}
[x_2,y_j] \in V_{2k_s+k_j},\\
x_{m} \in V_{mk_s}.
\end{cases}$$

Therefore, $2k_s \leq k_j \leq (n-p-2)k_s,$ which is only possible for $A=0,$ that is, $[x_2,y_j]=0, \ 1 \leq i \leq \frac{p}{2}.$

By applying the similar argumentations it can be proved that
$$\begin{array}{lll}
[x_i,y_j]=0=[y_j,x_i], &1 \leq i \leq n-p, & 1 \leq j \leq \frac{p}{2}, \qquad (i,j)\neq (1,j),\\{}
[x_j,z_i]=0=[z_i,x_j], & 1 \leq j \leq n-p,& 1 \leq i \leq \frac{p}{2},\\{}
[y_i,z_j]=0, &  1 \leq i ,j \leq \frac{p}{2}.&
\end{array}$$

Thanks to the Leibniz identity we get $[[z_i,y_j],x_1]=0.$ Since $[x_i,x_1]=0$ if and only if $x_i=x_{n-p}$, then
taking into account the product $[z_i,y_j]=(*)e_4+ \dots+ (*)e_{n-p},$ we obtain $[z_i,y_j]=A_{ij}x_{n-p}$  for some coefficients $A_{ij}$.

Furthermore we obtain $[z_i,y_i]=0, \ 1 \leq i \leq \frac{p}{2}$ by the properties of the gradation. Thanks to the Leibniz identity we derive $A_{ij}=A_{ji}, \ 1 \leq i,j \leq \frac{p}{2}.$


By following the same reasons we conclude $[y_i,y_j]=B_{ij}x_{n-p}$ with $B_{ii}=0.$


Finally, it is trivial to check that $[z_i,z_j]=0$ for $1 \leq i, j \leq \frac{p}{2},$ by using the Leibniz identity.

Summarizing, the law of the algebra is determined by the following products:
$$\ll: \begin{cases}
[x_i,x_1]=x_{i+1}, & 1 \leq i \leq n-p-1,\\
[x_1,y_i]=z_i, & 1 \leq i \leq \frac{p}{2},\\
[y_i,y_j]=B_{ij}x_{n-p}, & 1 \leq i,j \leq \frac{p}{2}, \  i \neq  j,\\
[z_i,y_j]=A_{ij}x_{n-p}, & 1 \leq i, j \leq \frac{p}{2}, \  i \neq j, \ A_{ij}=A_{ji},
\end{cases}$$
 with the conditions:
  \begin{equation} \label{parametros}
  \begin{cases}
  \hbox{if } B_{i_0j_0} \neq 0 \ \hbox{for \ some} \ i_0, j_0,  \ \hbox{then}  \ B_{i_0k}=0 \ \hbox{for all}\ \neq j_0 \ \hbox{and} \  B_{sj_0}=0 \ \hbox{for all}\ s \neq i_0,\\
   \hbox{if } A_{i_0j_0} \neq 0 \ \hbox{for \ some} \  i_0, j_0, \ \hbox{then} \ A_{i_0k}=0 \ \hbox{for all}\ k \neq j_0 \ \hbox{and} \ A_{sj_0}=0 \ \hbox{for all}\ s \neq i_0.
   \end{cases}
  \end{equation}

\

{\bf Case 1: $B_{ij}=0$ for all $i,j$.} Then the law of $\ll$ has the following form:
  $$\ll: \begin{cases}
[x_i,x_1]=x_{i+1}, & 1 \leq i \leq n-p-1,\\
[x_1,y_i]=z_i, & 1 \leq i \leq \frac{p}{2},\\
[z_i,y_j]=A_{ij}x_{n-p}, & 1 \leq i, j \leq \frac{p}{2}, \  i \neq j, \ A_{ij}=A_{ji},
\end{cases}$$
  where the parameters $A_{ij}$ satisfy the previous hypothesis.

  If all the parameters $A_{ij}$ are equal to zero, then it is easy to see that the algebra $\ll$ has maximum length. Indeed, putting
 $$V_i=\langle x_i\rangle, \ 1\leq i \leq n-p, \quad V_{n-p+2j-1}=\langle y_j \rangle, 1\leq j \leq \frac{p}{2}, \quad V_{n-p+2j}=\langle  z_j\rangle, 1 \leq j \leq \frac{p}{2},$$
we get
$$\ll=V_1\oplus V_2\oplus \dots \oplus V_{n-p}\oplus V_{n-p+1}\oplus V_{n-p+2}\oplus V_{n-p+3}\oplus V_{n-p+4}\oplus \dots \oplus  V_{n-1}\oplus V_{n}.$$ Thus we get the algebra $M^4(0).$

Let us assume now that there exists $A_{i_0j_0} \neq 0.$ Without loss of generality, we can suppose $A_{12}=0.$

Taking the following change of basis
$$y_i'=y_1-y_i, \quad y_{i+1}'=y_2+y_{i+1}, \quad z_i'=[x_1,y_i'], \quad z_{i+1}'=[x_1,y_{i+1}'], \ 3 \leq i \leq \frac{p}{2}-1,$$

we obtain the algebra, which is defined by the following products
   $$\ll: \begin{cases}
[x_i,x_1]=x_{i+1}, & 1 \leq i \leq n-p-1,\\
[x_1,y_i]=z_i, & 1 \leq i \leq \frac{p}{2},\\
[z_1,y_2]=x_{n-p},\\
[z_2,y_1]=x_{n-p}.
\end{cases}$$

Using the maximum length gradation properties and connectedness it is not difficult to check that this algebra admits the gradation of maximum length only for the cases of $k_s=\frac{n}{n-p}\in \mathbb{N}$ and $n$ even. The associated gradation is decomposed into direct sum of the following spaces:
$$\begin{array}{llll}
x_i \in V_{ik_s,} & &1 \leq i \leq n-p,&\\
y_1 \in V_{1},\qquad y_2\in V_{{(n-p-1)}k_s-1},& & &\\
z_1\in  V_{k_s+1}, \ \ z_2\in V_{{(n-p)k_s}-1},& & &\\
y_i\in V_{i-1}, & &3\leq i \leq k_s,&\\
z_i\in V_{k_s+i-1}, & &3 \leq i \leq k_s,&\\
y_{q(k_s-1)+i}\in V_{2qk_s-1+i},& & 1\leq q \leq \frac{n-p-4}{2}, & 2 \leq i \leq k_s,\\
z_{q(k_s-1)+i}\in V_{(2q+1)k_s-1+i},& & 1\leq q \leq \frac{n-p-4}{2}, & 2 \leq i \leq k_s,\\
y_{\frac{n-p-2}{2}(k_s-1)+i}\in V_{(n-p-2)k_s-1+i},& &2\leq i \leq k_s-1,&\\
z_{\frac{n-p-2}{2}(k_s-1)+i}\in V_{(n-p-1)k_s-1+i},& &2\leq i \leq k_s-1.&
\end{array}$$

\

{\bf Case 2: $\exists i_0, \ j_0$ such that $B_{i_0j_0} \neq 0$.} Making the basis transformation $y_1'=y_{i_0},$ $y_2'=y_{j_0},$ without loss of generality, one can assume that $B_{12}\neq 0.$ Further, applying the properties of gradation of maximum length, the conditions (\ref{parametros}) and the changes of basis, we arrive to the algebra of maximum length with the following table of multiplication:

$$\ll: \begin{cases}
[x_i,x_1]=x_{i+1}, & 1 \leq i \leq n-p-1,\\
[x_1,y_i]=z_i, & 1 \leq i \leq \frac{p}{2},\\
[y_1,y_2]=x_{n-p},
\end{cases}$$
whose associated graded spaces are
$$V_{-1}=\langle y_1 \rangle, \quad V_0=\langle z_1 \rangle, \quad V_i=\langle x_i \rangle, \quad 1 \leq i \leq n-p,$$
$$V_{n-p+2k+1}=\langle y_{k+2} \rangle, \ 0 \leq k \leq \frac{p}{2}-2, \quad V_{n-p+2k}=\langle z_{k+1}\rangle, \ 1 \leq k \leq \frac{p}{2}-1.$$

 \

 \

The description of the extension $\widetilde{M^2}$ is carried out in a similar way as for the extension $\widetilde{M^1}$.

\end{dem}

In the following theorem we prove that there is no algebra of maximum length among algebras from extension
$\widetilde{M^3}.$

\begin{thm}
Let $\ll$ be a $n$-dimensional $p$-filiform Leibniz algebra with $n-p \geq 4,$ $p \geq 4$ and $p$ odd. Then $\ll$ does not admit a gradation of maximum length.
\end{thm}

\begin{dem}
First of all, we shall denote $q=\lfloor\frac{p}{2}\rfloor$ for simplicity. Recall, the natural gradation of $M^3$ is $$L_1=\langle e_1,f_1,f_2, \dots, f_{q+1} \rangle \oplus L_2=\langle e_2,f_{q+2}, \dots , f_p \rangle \oplus L_{i}=\langle e_i \rangle, \ 3 \leq i \leq n-p.$$

By considering the law of $M^3$ we denote by new generators the following:
\begin{align*}
\widetilde{x_s}&=e_1+ \sum_{i=2}^{n-p} a_ie_i+\sum_{i=1}^p b_jf_j,\\
y_j&=f_j+\sum_{k=1}^{n-p}c_{kj}e_k+\sum_{k=1, k\neq j}^{\frac{p}{2}}d_{kj}f_k \quad \hbox{ for } 1 \leq j \leq q+1.
\end{align*}

Then we get
\small \begin{align*}
[\widetilde{x_s},\widetilde{x_s}]&=(1+b_{q+1})e_2+(*)e_3+\dots+(*)e_{n-p}+b_1f_{q+2}+ \dots + b_{q}f_p,\\
\underbrace{[[\widetilde{x_s},\widetilde{x_s}] \dots, \widetilde{x_s}]}_{i-times}&=(1+b_{q+1})^{i-1}e_i+(*)e_{i+1}+ \dots + (*)e_{n-p}, \quad 3 \leq i \leq n-p.\\
[\widetilde{x_s},y_1]&=(c_{11}+d_{q+1 1})e_2+(*)e_3+\dots + (*)e_{n-p}+ f_{q+2}+ d_{21}f_{q+3}+\dots+d_{q1}f_p,\\
[y_1,\widetilde{x_s}]&=c_{11}[\widetilde{x_s},\widetilde{x_s}],\\
[\widetilde{x_s},y_i]&=(c_{1i}+d_{q+1i})e_2+(*)e_3+\dots + (*)e_{n-p}+ d_{1i}f_{q+2}+\dots + f_{q+i+1}+\dots +d_{qi}f_p,\quad 2 \leq i \leq q,\\
[y_i,\widetilde{x_s}]&=c_{1i}[\widetilde{x_s},\widetilde{x_s}], \quad 2 \leq i \leq q+1,\\
[y_1,y_1]&=c_{11}[\widetilde{x_s},y_1],\\
[y_i,y_j]&=c_{1i}[\widetilde{x_s},y_j], \quad 1 \leq i,j \leq q+1, \ (i,j)\neq (1,1).
\end{align*}

 Without loss of generality, one can assume $1+b_{q+1} \neq 0$.
Then the homogeneous basis of $\ll$ is defined by the following vectors:
$$\begin{array}{ll}
x_1=\widetilde{x_s},&\\
x_i=[x_{i-1},\widetilde{x_s}], & 2 \leq i \leq n-p,\\
y_i,& 1 \leq i \leq q+1,\\
z_i=[x_1,y_i], &1 \leq i \leq q.
\end{array}$$

The matrix of the change of basis is as follows:
$$ \left( \begin{array}{ccccccccccccc}
1 & a_2 & a_3 & a_4 & \dots & a_{n-p} & b_1 & b_2 & \dots & b_{q+1} & b_{q+2} & \dots & b_p\\
0 & D_1 & (*) & (*) & \dots & (*) & 0 & 0 & \dots & 0 & b_1 & \dots & b_q\\
0 & 0 & D_2 & (*) & \dots & (*) & 0 & 0 & \dots & 0 & 0 & \dots & 0\\
0 & 0 & 0 & D_3 & \dots & (*)  & 0 & 0 & \dots & 0 & 0 & \dots & 0\\
0 & 0& 0 & 0 & \dots & D_{n-p-1} & 0 & 0 & \dots & 0 & 0 & \dots & 0\\
c_{11} & c_{21} & c_{31} & c_{41} & \dots & c_{n-p1} & 1 & d_{21} & \dots & d_{q+11} & d_{q+21} & \dots & d_{p1}\\
c_{12} & c_{22} & c_{32} & c_{42} & \dots & c_{n-p2} & d_{12} & 1 & \dots & d_{q+12} & d_{q+22} & \dots & d_{p2}\\
\vdots & \vdots &\vdots & \vdots &\dots & \vdots &\vdots & \vdots &\dots & \vdots &\vdots & \dots &\vdots\\
c_{1q+1} & c_{2q+1} & c_{3q+1} & c_{4q+1} & \dots & c_{n-pq+1} & d_{1q+1} & d_{2q+1} & \dots & 1 & d_{q+2q+1} & \dots & d_{pq+1}\\
0 & E_1 & (*) & (*)& \dots & (*) & 0 & 0 & \dots & 0 & 1 & \dots & d_{q1}\\
\vdots & \vdots &\vdots & \vdots &\dots & \vdots &\vdots & \vdots &\dots & \vdots &\vdots & \dots &\vdots\\
0 & E_q & (*) & (*)& \dots & (*) & 0 & 0 & \dots & 0 & d_{1q} & \dots & 1\\
\end{array} \right)$$
where $D_i=(1+b_{q+1})^i, \ 1 \leq i \leq n-p-1$ and $E_i=c_{1i}+d_{q+1i}, \ 1 \leq i \leq q.$

The gradation associated to the above basis is
$$V_{k_s}\oplus V_{2k_s}\oplus \dots \oplus V_{(n-p)k_s}\oplus V_{k_1}\oplus \dots \oplus V_{k_{q+1}}\oplus V_{k_1+k_s} \oplus \dots \oplus V_{k_{q}+k_s}.$$

Let us assume that this gradation has maximum length.

Consider
$$[x_2,y_i]=(1+b_{q+1})(c_{1i}+d_{q+1i})e_3+(*)e_4+ \dots+(*)e_{n-p},\qquad \ 1\leq i\leq q.$$

Therefore we conclude that $[[x_1,x_1],y_i]=Ax_3,$ with $A \in \mathbb{C}.$

On the other hand, by considering the properties of the gradation we have
$$\begin{cases}
[x_2,y_i] \in V_{2k_s+k_i}, \hbox{ for } 1 \leq i \leq q,\\
x_3 \in V_{3k_s},
\end{cases}$$
that is, we have either $A=0$ or $k_s=k_i.$ The last equality contradicts the assumption of maximum length, thus $A=0,$ i.e, we obtain
\begin{equation}\label{eq1}
c_{1i}+d_{q+1i}=0, \quad 1 \leq i \leq q.
\end{equation}

From the product $$[x_2, y_{q+1}]=(1+b_{q+1})(c_{1q+1}+1)e_3+(*)e_4+\dots+(*)e_{n-p},$$ we conclude $c_{1q+1}=-1.$

Finally, it contradicts the assumption of maximum length by comparison with the following equalities:
$$[y_{q+1},x_1]=-(1+b_{q+1})e_2+(*)e_3+\dot+(*)e_{n-p}-b_1f_{q+2}- \dots-b_qf_p=-x_2.$$

By means of the gradation, it leads to $V_{k_s+k_{q+1}}=V_{2k_s},$ i.e, $k_s=k_{q+1},$ which contradicts the properties (\ref{grad}).

  \end{dem}



\begin{thebibliography}{12}
\bibitem{Jobir} {\rm  Adashev J.Q., Ca\~{n}ete E.M.,}{\it Derivations of the quasi-filiform Leibniz algebras of maximum length},  Uzbek Math. J.,  2, 2011, 3--14.

\bibitem{Omirov1} {\rm Ayupov Sh.A., Omirov B.A.,} {\it On some classes of nilpotent
Leibniz algebras}, (Russian) Sib. Mat. Zh., 42(1), 2001,
18--29; translation in  Sib. Math. J., 42(1), 2001,
15--24.


\bibitem{Bar} {\rm Barnes D.W.}, {\it On {L}evi's theorem for {L}eibniz algebras}, Bull.
  Aust. Math. Soc., 86(2), 2012, 184--185.

\bibitem{Pastor}{\rm Cabezas J.M., Pastor E.,} {\it Naturally graded $p$-filiform Lie algebras in arbitrary finite dimension}, J. Lie Theory, 15, 2005, 379--391.

\bibitem{NGMio} {\rm Camacho L.M., Ca\~{n}ete E.M., G\'{o}mez J.R., Omirov B.A.,} {\it Quasi-filiform Leibniz algebras of maximum length,} Sib. Math. J., 52(5), 2011, 840--853.

\bibitem{3-filiform}{\rm Camacho L. M., Ca\~{n}ete E.M, G\'{o}mez J.R., Omirov B.A.,}
 {\it 3-filiform Leibniz algebras of maximum length, whose naturally graded algebras are Lie algebras}, Linear Multilinear Alg., 59(9), 2011, 1039--1058.

\bibitem{3-filiformLeib}{\rm Camacho L. M., Ca\~{n}ete E.M, G\'{o}mez J.R., Omirov B.A.,}
{\it 3-filiform Leibniz algebras of maximum length.}, Sumitted to Journal of Algebra (2013), arXiv:1310.6539v1.


\bibitem{NGp-F} {\rm Camacho L.M., G\'{o}mez J.R., Gonz\'{a}lez A.J., Omirov B.A.} {\it The classification of naturally graded $p$-filiform Leibniz algebras,} Commun. Alg., 39(1), 2011, 153--163.

\bibitem{mitesis}{\rm Ca\~{n}ete E.M.} {\rm Algebras de Leibniz de longitud maxima,} PhD Thesis, Universidad de Sevilla, 2012. (http://www.educacion.es/teseo).




\bibitem{Casas}{\rm Casas J.M., Ladra M., Omirov B.A., Karimjanov I.A.} {\it Classification of solvable Leibniz algebras with null-filiform nilradical}, Linear Multilinear Alg., 61(6), 2012, 758--774.

\bibitem{Jimenez1}{\rm G\'{o}mez J.R., Jim\'{e}nez-Merch\'{a}n A., Reyes J.} {\it Maximum length filiform Lie algebras}, Extracta Mathematicae, 16(3), 2001, 405--421.


\bibitem{Gomez}{\rm G\'{o}mez J.R., Jim\'{e}nez-Merch\'{a}n A., Reyes J.}, {\it Quasi-filiform Lie algebras of maximum length,} Linear Algebra and its Applications, 335, 2001, 119--135.

\bibitem{Jacobson}{\rm Jacobson N.,} {\it Lie algebras,} Interscience Tracts in Pure and Applied Mathematics, No. 10, Interscience Publishers (a division of John Wiley and Sons), New York-London, 1962.


\bibitem{loday} {\rm Loday J.L.,} {\it Une version non commutative des
alg$\acute{e}$bres de Lie: les alg$\acute{e}$bres de Leibniz.}
Ens. Math., 39, 1993, 269--293.

\bibitem{Mal'cev}{\rm Mal'cev A.I.,} {\it Solvable Lie algebras,} Amer. Math. Soc. Translation, 1950, p. 36.

\bibitem{Mub} {\rm Mubarakzjanov G.M.} {\it On solvable {L}ie algebras,} (Russian), Izv. Vys\v s. U\v cehn.
  Zaved. Matematika 1963, 114--123.



\end{thebibliography}
\end{document}